\documentclass[12pt,reqno]{amsart}
\usepackage{graphicx}
\usepackage{amsmath}
\usepackage{amssymb}
\usepackage{amsfonts}
\usepackage{euscript}
\usepackage[latin1]{inputenc}
\usepackage[T1]{fontenc}
\numberwithin{equation}{section}
\makeatletter\@addtoreset{equation}{section}

\newcommand{\norm}[1]{{\left\|{#1}\right\|}}
\newcommand{\set}[1]{{\left\{{#1}\right\}}}
\newcommand{\scal}[1]{{\left\langle{#1}\right\rangle}}
\newcommand{\C}{\mathbb C}
\newcommand{\R}{\mathbb R}

\newcommand{\fin}{\hfill  \rule {2mm}{1.5mm}}
\newcommand{\Pag}{\frac{\partial}{\partial \bg}}

\newcommand{\mi}{m\curlyvee n}
\newcommand{\ma}{m\curlywedge n}
\newcommand{\Pa}{\frac{\partial}{\partial {z}}}
\newcommand{\bPa}{\frac{\partial}{\partial {z^*}}}
\newcommand{\Pan}{\frac{\partial^n}{\partial {z^n}}}
\newcommand{\bPam}{\frac{\partial^m}{\partial {{z^*}^m}}}
\newcommand{\PbP}{\frac{\partial^2}{\partial z\partial z^*}}
\newcommand{\PnPm}{\frac{\partial^{m+n}}{\partial z^n\partial {z^*}^m}}
\newcommand{\bz}{z^*}
\newcommand{\g}{\xi}
\newcommand{\bg}{\xi^*}
\newcommand{\fa}{\mathfrak{A}_{\nu,\xi}}
\newcommand{\bfa}{\mathfrak{A}_{\nu,\xi}^*}
\newcommand{\La}{\mathfrak{L}_{\nu,\xi}}
\newtheorem {theorem}{Theorem}[section]
\newtheorem {definition}[theorem]{Definition}
\newtheorem {lemma}[theorem]{Lemma}
\newtheorem {proposition}[theorem]{Proposition}
\newtheorem {remark}[theorem]{Remark}

\begin{document}

\title{A class of generalized  complex Hermite polynomials}

\author{Allal Ghanmi}

\address{{\sf Allal Ghanmi ({\it Current address})}
\newline Center for Advanced Mathematical Sciences,
P. O. Box 11-0236,
\newline College Hall, 4th Floor,
American University of Beirut,
Beirut, Lebanon
\newline {\it \hspace*{.6cm} Permanent address}
\newline Secteur E, N 380, Hay Errahma, 11 000 Sal\'e, Maroc
}
 \email{allalghanmi@gmail.com}

 \thanks{The author would like to acknowledge the financial support
 of the Arab Regional Research Program  2006-2007 via CAMS, AUB.}

\date{August 20, 2007}
\maketitle

\begin{abstract}
A class of generalized complex polynomials of Hermite type,
suggested by a special magnetic Schr\"odinger operator, is
introduced and some related basic properties  are discussed.
\end{abstract}

\section{ Introduction}
Let $\C$ be the space of complex numbers $z=x+iy$; $x,y\in \R$, and
denote by $\partial/{\partial {z}}$ and $\partial/{\partial {z^*}}$
the derivation with respect to the variable $z$ and its conjugate
$\bz$, respectively; i.e.,
\[
\Pa=\frac 12 \Big(\frac\partial{\partial x}-i\frac\partial{\partial
y}\Big) \qquad \mbox{and} \qquad \bPa=\frac 12
\Big(\frac\partial{\partial x}+i\frac\partial{\partial y}\Big).\]
For given fixed $\nu>0$ and $\g\in \C$, we set
$S_{\nu,\g}=S_{\nu,\g}(z)= \nu z+\g$ and then consider the elliptic
selfadjoint second order differential operator $\La$ given
explicitly in complex coordinate $z$ by
      \begin{equation}\label{AutomLap1}
      \La := -\frac 14 \set{4\PbP +2(S_{\nu,\g}\Pa
      -S_{\nu,\g}^*\bPa ) -|S_{\nu,\g}|^2}.
      \end{equation}
Note that for $\g=0$, it gives rise to the special Hermite operator
$\mathfrak{L}_{\nu}$ (called also twisted Laplacian),
      \begin{equation}\label{LanHam}
      \mathfrak{L}_{\nu}=-\frac 14 \set{4\PbP +2\nu(z\Pa -\bz \bPa)-\nu^2|z|^2},\nonumber
      \end{equation}
which describes in physics a nonrelativistic quantum particle moving
on the plane under the action of an external constant magnetic field
applied perpendicularly. The associated eigenfunctions are known to
be expressible in terms of the complex Hermite polynomials
\cite{Sh87,Thangavelu93,IntInt06},
      \begin{equation}
      \label{CompHerPoly} H^{m,n}(z,\bz):= (-1)^{m+n} e^{|z|^2}\PnPm
      e^{-|z|^2} .
      \end{equation}
Such polynomials form a complete orthogonal system of the Hilbert
space $L^2(\C;e^{-\nu|z|^2}d\lambda)$, where $d\lambda$ being the
Lebesgue measure on $\C$, and appear as an essential tool in many
area of mathematics and physics. They have been studied by many
authors, notably by Shigekawa \cite{Sh87}, Thangavelu
\cite{Thangavelu93}, W\"{u}nsche \cite{Wunsche99-1,Wunsche99-2},
Dattoli \cite{Dattoli03} and more recently by Intissar and Intissar
\cite{IntInt06}.\\
 Our aim in the present paper is to discuss some
basic properties of a general class of complex polynomials
$\mathfrak{G}_{\nu}^{m,n}(z,\bz|\g)$ of Hermite type suggested by
the Laplacian $\La$ such that the associated functions
      \[
      \mathfrak{g}_{\nu}^{m,n}(z,\bz|\g)=    e^{-\frac 12 \bz S_{\nu,\g}}
      \mathfrak{G}_{\nu}^{m,n}(z,\bz|\g),
      \]
are solutions of the eigenvalue problem $\La \psi = \mu \psi$, for
$\psi\in  \mathcal{C}^\infty(\C)$ and $\mu \in \C$. More precisely,
we show that the involved polynomials satisfy the Rodriguez type
formula (\ref{RodriguezTypeFormula}) and can be expressed as
binomial sum of the complex Hermite polynomials (\ref{CompHerPoly}),
see identity (\ref{Identity21}). Other properties of these
polynomials such as three-term recursion relations and differential
equations which they obey are obtained in Section 3. Also the
generating function as well as the explicit series expansion and
representation by the means of the confluent hypergeometric
function, Laguerre and complex Hermite polynomials are derived
(Section 4). Furthermore, the weak orthogonal property of these
polynomials is discussed and their norms are explicitly determined
(Section 5). We conclude by giving some new identities for the usual
complex Hermite and Laguerre polynomials and by illustrating some
obtained results making use of matrix representation (Section 6).

\section{ A Rodriguez type formula for  $\mathfrak{G}_{\nu}^{m,n}(z,\bz|\g)$}
In this section we introduce the polynomials
$\mathfrak{G}_{\nu}^{m,n}(z,\bz|\g)$ using a classical approach (see
for instance \cite{Thangavelu93}). For fixed  $\nu>0$ and $\g\in\C$,
we denote by $\fa$ and $ \bfa $ the first order differential
operators given respectively by
      \[
      \fa   = \bPa +\frac 12 S_{\nu,\g} \qquad \mbox{and}  \qquad \bfa =
      -\Pa +\frac 12 S_{\nu,\g}^{*},
      \]
where $S_{\nu,\g}(z)=\nu z +  \g$. Then, one check easily the
following algebraic relationships
      \begin{align}
      \fa \bfa  = \La + \frac 12\nu , \qquad
      \bfa \fa  = \La - \frac 12\nu.\label{AlgId2}
      \end{align}
Hence, it can be shown that the null space, $\ker(\fa )= \set{\psi \in \mathcal{C}^\infty(\C) ; ~\fa \psi
= 0 },$ of the operator $\fa $ coincides with the eigenspace
\[\mathcal{E}_{0}(\La):=\set{\psi \in \mathcal{C}^\infty(\C);\quad
\La\psi = \frac\nu 2 \psi},\] and therefore the functions
$\psi_{\nu,\g}^m(z,\bz)$ given  explicitly by
      \begin{equation}\label{LikeHermiteFun0}
      \psi_{\nu,\g}^m(z,\bz) := z^m e^{-\frac 12 \bz S_{\nu,\g}}= z^m e^{-\frac 12( \nu|z|^2 + \g
      \bz)}; \quad m=0,1,2, \cdots,
      \end{equation}
span linearly $\mathcal{E}_{0}(\La)$. Moreover, the functions
      \begin{equation}\label{LikeHermiteFun1}
      \mathfrak{g}_{\nu}^{m,n}(z,\bz|\g) := \Big[[\bfa ]^n\psi_{\nu,\g}^m\Big](z,\bz),
      \end{equation}
are eigenfunctions of  $\La$ with (the Landau levels) $\nu(n+\frac
12)$; $n=0,1,2, \cdots $, as corresponding eigenvalues. The involved
operator $[\bfa ]^n $ is given by
      \begin{align*}
      [\bfa ]^n \varphi  = (-1)^n e^{\frac 12 z S_{\nu,\g}^*}\Pan \Big[
      e^{-\frac 12 z S_{\nu,\g}^*}\varphi\Big].
      \end{align*}
Hence, one can rewrite the functions
$\mathfrak{g}_{\nu}^{m,n}(z,\bz|\g)$ as follows
      \begin{align}
      \mathfrak{g}_{\nu}^{m,n}(z,\bz|\g)
      &= (-1)^n  e^{-\frac 12 \bz S_{\nu,\g}} e^{\nu|z|^2+\Re e\scal{z,\g}}\Pan
      \Big(z^m e^{-\nu|z|^2-\Re e\scal{z,\g}}\Big)\label{LikeHermiteFunA1}\\
      &= (-1)^{n}e^{-\frac 12\bz S_{\nu,\g}}e^{\nu|z|^2 +\frac \bg 2z}\Pan
      \Big(z^m e^{-\nu|z|^2-\frac \bg 2z}\Big)\label{LikeHermiteFunA2}
      \end{align}
and introduce   $\mathfrak{G}_{\nu}^{m,n}(z,\bz|\g)$ to be the
polynomial of degree $m$ in $z$ and degree $n$ in  $\bz$ defined by
      \begin{align}
      \mathfrak{G}_{\nu}^{m,n}(z,\bz|\g)
      & =   e^{\frac 12 \bz S_{\nu,\g}}\mathfrak{g}_{\nu}^{m,n}(z,\bz|\g) \nonumber \\
      &= (-1)^{n} e^{\nu|z|^2 +\frac \bg 2z}\Pan \Big(z^m e^{-\nu|z|^2-\frac \bg
      2z}\Big).\label{LikeHermitePolyC}
      \end{align}
Thus, we have

\begin{proposition}
The polynomials $\mathfrak{G}^{m,n}(z,\bz)$ satisfy the following
Rodriguez type formula
      \begin{equation}\label{RodriguezTypeFormula}
      \mathfrak{G}_{\nu}^{m,n}(z,\bz|\g) = \frac{(-1)^{m+n}}{\nu^m} e^{\nu|z|^2+\frac \bg 2 z}
      \PnPm\Big(e^{-\nu|z|^2 -\frac \bg 2z}\Big).
      \end{equation}
\end{proposition}

\noindent The above expression can serves as definition for this
class of complex polynomials of Hermite type.

\begin{definition}
We call $\mathfrak{G}_{\nu}^{m,n}(z,\bz|\g)$ redefined by
(\ref{RodriguezTypeFormula}) generalized complex Hermite polynomials
(GCHP).
\end{definition}

\begin{remark}
For $\g=0$ and $\nu=1$, the polynomials
$\mathfrak{G}_{1}^{m,n}(z,\bz|0)$ reduce further to the complex
Hermite polynomials $H^{m,n}(z,\bz)$ as given by
(\ref{CompHerPoly}).
\end{remark}

\section{ Related basic properties}
In this section we investigate some properties of the polynomials
$\mathfrak{G}_{\nu}^{m,n}(z,\bz|\g)$ introduced above. We begin with
the following
\begin{proposition}
The polynomials
$\mathfrak{G}_{\nu}^{m,n}(z,\bz|\g)$ satisfy the following identity
      \begin{align}
      \mathfrak{G}_{\nu}^{m,n}(z,\bz|\g)
      = \frac{n!}{\sqrt{\nu}^m}\sum_{j=0}^n
      \frac{\sqrt{\nu}^j}{j!}
      \frac{(\bg/2)^{n-j}}{(n-j)!}H^{m,j}(\sqrt{\nu}
      z,\sqrt{\nu}\bz).
      \label{Identity21}
      \end{align}
\end{proposition}

\noindent {\sf Proof.} The proof of (\ref{Identity21}) relies
essentially on (\ref{RodriguezTypeFormula}). Indeed, it follows by
rewriting the polynomial $\mathfrak{G}_{\nu}^{m,n}(z,\bz|\g)$ as
      \begin{align*}
      \mathfrak{G}_{\nu}^{m,n}(z,\bz|\g) &= \frac{(-1)^{m}}{\nu^m}
      e^{\nu|z|^2+\frac \bg 2 z} \bPam\Big(\big(\nu\bz +\frac \bg 2\big)^n
      e^{-\nu|z|^2 -\frac \bg 2 z}\Big)\\
      &= (-1)^{m} e^{\nu|z|^2}
      \bPam\Big(\big(\bz +\frac \bg {2\nu}\big)^n e^{-\nu|z|^2}\Big)
      \end{align*}
 and next by making use of the binomial formula to expand $\big(\bz +\frac \bg
 {2\nu}\big)^n$.

 \fin

Furthermore, we have
\begin{proposition}\label{Prop3TRR}
 The polynomials $\mathfrak{G}_{\nu}^{m,n}(z,\bz|\g)$ satisfy
 the following three-terms recursion relations
      \begin{align}
      & \mathfrak{G}_{\nu}^{m,n+1}(z,\bz|\g) = -  m
      \mathfrak{G}_{\nu}^{m-1,n}(z,\bz|\g)+ \Big(\nu\bz +  \frac \bg
      2       \Big)\mathfrak{G}_{\nu}^{m,n}(z,\bz|\g), \label{RecursionRelation}\\
      & \mathfrak{G}_{\nu}^{m+1,n}(z,\bz|\g) = -  n
      \mathfrak{G}_{\nu}^{m,n-1}(z,\bz|\g)+
      z\mathfrak{G}_{\nu}^{m,n}(z,\bz|\g). \label{RecursionRelation2}
      \end{align}
\end{proposition}

\noindent {\sf Proof.} By applying the operator $(-1)^{n+1} e^{
\nu|z|^2 + \frac\bg 2 z }\partial^n/\partial z^n $ to the both sides
of the following elementary fact
      \begin{equation*}
      \label{ElemFact}
      \Pa ( z^m e^{a|z|^2+bz} ) = m z^{m-1}e^{a|z|^2+bz} +(a \bz + b) z^me^{a|z|^2+bz}
      \end{equation*}
with $a=-\nu, b=-\g/2$, and the use of (\ref{RodriguezTypeFormula})
one obtains (\ref{RecursionRelation}). \\
Equation (\ref{RecursionRelation2}) holds by equating the right hand
sides of
      \[
      2\Pag \mathfrak{G}_{\nu}^{m,n}(z,\bz|\g) = z\mathfrak{G}_{\nu}^{m,n}(z,\bz|\g)
      - \mathfrak{G}_{\nu}^{m+1,n}(z,\bz|\g).    \eqno{(a)\qquad}
      \]
and
      \[
      2\Pag \mathfrak{G}_{\nu}^{m,n}(z,\bz|\g)
      =n\mathfrak{G}_{\nu}^{m,n-1}(z,\bz|\g) .   \eqno{(b)\qquad}
      \]
The fact (a) is obtained by differentiating both sides of the
Rodriguez type formula (\ref{RodriguezTypeFormula}) with respect to
$\bg$. While (b) can be handled from different ways, particularly by
the use of the identity (\ref{Identity21}) or also from the
generating function (\ref{GeneratingFunction}) given below. \fin

\begin{remark}\label{Rem3TRR}
i) Combination of (\ref{RecursionRelation}) and
(\ref{RecursionRelation2}) yields
      \begin{equation}
      \label{RecursionRelation3} {{     (m-n)
      \mathfrak{G}_{\nu}^{m,n}(z,\bz|\g) =
      -z\mathfrak{G}_{\nu}^{m,n+1}(z,\bz|\g)
      + (\nu\bz+\frac{\bg}2)\mathfrak{G}_{\nu}^{m+1,n}(z,\bz|\g).}}
      \end{equation}
From which we deduce
      \begin{equation} \label{RecursionRelation4}
      {{z\mathfrak{G}_{\nu}^{m,m+1}(z,\bz|\g)
      =(\nu\bz+\frac{\bg}2)\mathfrak{G}_{\nu}^{m+1,m}(z,\bz|\g).}}
      \end{equation}

ii) By taking $\g=0$ and $\nu=1$ in the previous obtained recursion
relations (\ref{RecursionRelation})-(\ref{RecursionRelation4}), one
deduce the following ones for the usual complex Hermite polynomials
\cite[Eqns in (2.9)]{Wunsche99-2},
      \begin{align*}
      & H^{m,n+1}(z,\bz) =  \bz H^{m,n}(z,\bz)-m H^{m-1,n}(z,\bz)\\
      & H^{m+1,n}(z,\bz)= z H^{m,n}(z,\bz) - n H^{m,n-1}(z,\bz)\\
      & (m-n) H^{m,n}(z,\bz) =-zH^{m,n+1}(z,\bz) + \bz H^{m+1,n}(z,\bz)\\
      & zH^{m,m+1}(z,\bz) =\bz H^{m+1,m}(z,\bz).
      \end{align*}
\end{remark}

Now, using (\ref{RecursionRelation}) and (\ref{RecursionRelation2}),
one can show  that the polynomials
$\mathfrak{G}_{\nu}^{m,n}(z,\bz|\g)$ obey second and first order
differential equations. Namely, we have

\begin{proposition}
The polynomials $\mathfrak{G}_{\nu}^{m,n}(z,\bz|\g)$ are solutions
of the following second order differential equations
      \begin{equation}
      \label{DiffEq}
      - \PbP +\nu z \Pa =\nu m \quad \mbox{and} \quad
      - \PbP + (\nu\bz+\frac{\bg}2)\bPa = \nu n.
      \end{equation}
Therefore, they satisfy the first order differential equation
      \begin{equation}\label{DiffEq3}
      \nu z \Pa - (\nu\bz+\frac{\bg}2)\bPa = \nu (m-n).
      \end{equation}
\end{proposition}

\noindent {\sf Proof.} To get the first equation in (\ref{DiffEq}),
we differentiate (\ref{RecursionRelation2}) w.r.t. the variable $z$
and next use the well established facts
      \[
      \Pa \mathfrak{G}_{\nu}^{m+1,n}(z,\bz|\g)= (m+1)
      \mathfrak{G}_{\nu}^{m,n}(z,\bz|\g) \eqno{(c)\qquad}
      \]
and
      \[
      n\nu \mathfrak{G}_{\nu}^{m,n-1}(z,\bz|\g)=\bPa
      \mathfrak{G}_{\nu}^{m,n}(z,\bz|\g).  \eqno{(d)\qquad}
      \]
The second equation in (\ref{DiffEq}) can be handled in a similar
way making use of (\ref{RecursionRelation}). It can also be obtained
from the fact that the functions
      \begin{align}\label{EigenF}
      \mathfrak{g}_{\nu}^{m,n}(z,\bz|\g)= e^{-\frac 12 \bz S_{\nu,\g}}
      \mathfrak{G}_{\nu}^{m,n}(z,\bz|\g) ; \qquad m=0,1,2, \cdots,
      \end{align}
are eigenfunctions of the magnetic Schr\"odinger operator $\La$
given by (\ref{AutomLap1}) with $\nu(n+\frac  12)$ as associated
eigenvalues. \fin

\begin{remark}  From Eqns (b) and (d), we deduce that
      \begin{align}\label{Deribzzeta}
      {{\bPa \mathfrak{G}_{\nu}^{m,n}(z,\bz|\g) =
      2\nu\Pag\mathfrak{G}_{\nu}^{m,n}(z,\bz|\g)}}.
      \end{align}
\end{remark}

\section{Expansion series and relationship to some special functions}
 We begin by  realizing the
polynomials $\mathfrak{G}_{\nu}^{m,n}(z,\bz|\g)$ as iteration of the
monomial $z^m$ via a first order differential operator. Precisely,
we have

\begin{proposition}
The polynomials
$\mathfrak{G}_{\nu}^{m,n}(z,\bz|\g)$ are the excited states of $z^m$
by the first order differential operator $-\Pa +\nu \bz +\frac \bg
2$. Precisely
      \begin{equation}
      \bigg(-\Pa +\nu \bz +\frac \bg 2\bigg)^n (z^m)
      =\mathfrak{G}_{\nu}^{m,n}(z,\bz|\g). \label{It4}
      \end{equation}
\end{proposition}

\noindent {\sf Proof.} The  above assertion follows by successive
application of the following fact
      \begin{equation}
      \bigg(-\Pa +\nu \bz +\frac \bg 2\bigg)\mathfrak{G}_{\nu}^{j,k}(z,\bz|\g)=
      \mathfrak{G}_{\nu}^{j,k+1}(z,\bz|\g)     \label{It5}
      \end{equation}
with $j=m$ and $k=0$, keeping in mind that
$\mathfrak{G}_{\nu}^{m,0}(z,\bz|\g)=z^m$. Note here that (\ref{It5})
is in fact equivalent to the recursion relation
(\ref{RecursionRelation}) thanks to Eqn. (c).\fin

\begin{remark}
 The result (\ref{It4}) is similar to the one
obtained in \cite[Eqn. (2.2)]{IntInt06} for the complex Hermite
polynomials $H^{m,n}(z,\bz)$.
\end{remark}

Therefore, for every positive integers $m,n$, we deduce  easily from
(\ref{It4}) (or also (\ref{RodriguezTypeFormula})) that
      \begin{align}
      & \mathfrak{G}_{\nu}^{m,0}(z,\bz|\g) =  z^m  \label{Example11} \\
      & \mathfrak{G}_{\nu}^{0,n}(z,\bz|\g) =\Big(\nu\bz + \frac \bg 2 \Big)^n. \label{Example12}
      \end{align}
Also for every given integers $m\geq 1$ and $n\geq 1$, one obtain
      \begin{align}
      & \mathfrak{G}_{\nu}^{m,1}(z,\bz|\g)= z^{m-1}\Big[z\Big(\nu\bz + \frac \bg 2\Big) -m \Big]
      \label{Example21}\\
      & \mathfrak{G}_{\nu}^{1,n}(z,\bz|\g) = \Big(\nu\bz +\frac \bg 2 \Big)^{n-1} \Big[ z\Big(\nu\bz + \frac \bg 2 \Big) - n\Big].
      \label{Example22}
      \end{align}
 More generally, the expansion series of the GCHP is given by

\begin{proposition}
Denote by $\mi$  the minimum of $m$ and $n$. Then
      \begin{align}
      \mathfrak{G}_{\nu}^{m,n}(z,\bz|\g)
      &= m!n!     \sum_{j=0}^{\mi }  \frac{(-1)^j}{j!}\frac{z^{m-j}}{(m-j)!}
      \frac{(\nu \bz +\frac \bg 2)^{n-j}}{(n-j)!}                                                     \label{LikeHermitePoly21} \\
      &= m!n! \nu^n \sum_{j=0}^{\mi }   \sum_{k=0}^{n-j} \frac{(-1/\nu)^j}{j!}
      \frac{z^{m-j}}{(m-j)!} \frac{{\bz}^{k}}{k!} \frac{( \bg /{2\nu})^{n-j-k}}{(n-j-k)!}
      \label{LikeHermitePoly24}
      \end{align}
which we rewrite also as
      \begin{align}
      \mathfrak{G}^{m,n}_\nu(z,\bz|\g)
      =  &(-1)^{m} m!n!\frac{\nu^n}{\nu^m} \times \label{LikeHermitePoly25}\\
      &\times\sum_{l=m-(\mi) }^{m}\sum_{k=0}^{l+n-m} \frac{(-\nu)^l}{(m-l)!}
      \frac{z^l}{l!} \frac{{\bz}^{k}}{k!} \frac{(\bg /{2\nu})^{n-m+l-k}}{(n-m+l-k)!}    \nonumber  .
      \end{align}
\end{proposition}

\noindent {\sf Proof.} Such expansion series can be obtained by
direct computation using (\ref{It4}) or also
(\ref{LikeHermitePolyC}) together with the application of the
Leibnitz formula for the n$^{th}$ derivative of a product. \fin

\begin{remark}
The monomial of the lowest degree in the expansion above of
$\mathfrak{G}_{\nu}^{m,n}(z,\bz|\g)$ is
      \[
      (-1)^{\mi }\frac{(\ma)!}{|m-n|!} \Big(\frac \bg 2\Big)^{n-(\mi)} z^{m-(\mi)},
      \]
where $\ma$ denotes the maximum of $m$ and $n$, so that the analogue
of the Hermite numbers for these polynomials are 
      \begin{align}
      \mathfrak{G}_{\nu}^{m,n}(0,0|\g) =  \left\{\begin{array}{lll}
                       0         & \quad \mbox{if} \quad  m > n \\
                       (-1)^m m! & \quad \mbox{if} \quad  m = n \\
                       (-1)^m n!\frac{(\bg /2)^{n-m}}{(n-m)!} & \quad \mbox{if} \quad  n>m
      \end{array}\right. .\label{HermiteNum}
      \end{align}
\end{remark}

We conclude this section by pointing out, from
(\ref{LikeHermitePoly21}), that the polynomials
$\mathfrak{G}_{\nu}^{m,n}(z,\bz|\g)$ are linked to the complex
Hermite polynomials by
      \begin{align}\label{LinkDattoli1}
      \mathfrak{G}_{\nu}^{m,n}(z,\bz|\g) = H^{m,n}(z,\nu\bz+\frac \bg 2)
      = h_{m,n}(z,\nu\bz+\frac \bg 2|-1),
      \end{align}
where the notation $h_{m,n}(z,\bz|\tau)$ is used by Dattoli in
\cite{Dattoli03} to mean
\[h_{m,n}(z,\bz|\tau) := m!n!\sum_{j=0}^{\mi } \frac{\tau^j}{j!}
\frac{z^{m-j}}{(m-j)!} \frac{{\bz}^{n-j}}{(n-j)!}. \] \noindent We
should note here that the most established properties of the
polynomials $\mathfrak{G}_{\nu}^{m,n}(z,\bz|\g)$ can be deduced
formally by replacing $\bz$ by $\nu\bz+\frac \bg 2$ in the different
known properties of $H^{m,n}(z,\bz)=h_{m,n}(z,\bz|-1)$. But in
general there is no reason to the obtained results be a strict
consequence of the prescription $\bz
\leftrightsquigarrow\nu\bz+\frac \bg 2$. Nevertheless, one can
deduce from \cite[Eq. (31)]{Dattoli03} the following generating
function, which can be handled directly using the Taylor expansion
series of the function on the right hand side.

\begin{proposition}
We have
      \begin{equation}\label{GeneratingFunction}
      \sum_{m=0}^\infty\sum_{n=0}^\infty \frac{u^m}{m!}\frac{v^n}{n!}
      \mathfrak{G}_{\nu}^{m,n}(z,\bz|\g) =e^{u z +v(\nu\bz+\frac \bg 2) - uv}
      \end{equation}
\end{proposition}

\noindent Furthermore, using the dependence (\ref{LinkDattoli1}) (or
Eqn. (\ref{It4})) together with the fact \cite[Eqn.
(2.3)]{IntInt06}, one can check that the
$\mathfrak{G}_{\nu}^{m,n}(z,\bz|\g)$  can be written also in terms
of the confluent hypergeometric function ${_1F_1}(a;c;x)$,
      \[
      {_1F_1}(a;c;x)= \frac{\Gamma(c)}{\Gamma(a)}\sum_{k=0}^\infty \frac{\Gamma(a+j)}{\Gamma(a+j)}\frac{x^j}{j!}
      \]
as follows

\begin{proposition}
We have
      \begin{align}
      \mathfrak{G}_{\nu}^{m,n}(z,\bz|\g)
      &= \frac{(-1)^{\mi}(\mi)!}{|m-n|!} z^{(\ma) - n}
      (\nu\bz +\frac{\bg}{2})^{(\ma) -m}\times \nonumber \\
      & \quad  \times{_1F_1}(-(\mi) ; |m-n|+1 ;
      \nu|z|^2+\frac{\bg}{2}z). \label{LikeHyper}
      \end{align}
\end{proposition}

\begin{remark}
Using  the known fact
\[ {_1F_1}(-s ; c +1 ;x) = \frac{\Gamma(s+1)\Gamma(c+1)}{\Gamma(s+c +1)}L^{c}_s(x): \quad s=0,1,2, \cdots, \]
one can express the polynomials $\mathfrak{G}_{\nu}^{m,n}(z,\bz|\g)$
in terms of the Laguerre polynomials $L^{c}_s(x)$ as
      \begin{align}
      \mathfrak{G}_{\nu}^{m,n}(z,\bz|\g)
      = (-1)^{\mi}(\mi)!  & z^{(\ma) -n} (\nu\bz +\frac{\bg}{2})^{(\ma)
      -m}\times \nonumber\\
      & \times L_{\mi}^{|m-n|}(\nu|z|^2+\frac{\bg}{2}z).\label{LikeLaguerre}
      \end{align}
\end{remark}

\section{Weak orthogonality.}
Denote by $d\lambda$ the Lebesgue measure on $\C$. Let $\omega$ be a
positive weighting function on $\C$ and define $\scal{f , g
}_{\omega}$ by
      \begin{align}\label{WeakOrtho2}
      \scal{f  , g}_{\omega}:= \int_{\C}f(z)[g(z)]^{*}
      \omega(z)d\lambda(z).
      \end{align}
Then, we state

\begin{proposition}\label{PropId1}
Suppose that the system
$\set{\mathfrak{G}^{m,n}:=\mathfrak{G}_{\nu}^{m,n}(\cdot,\cdot|\g)}_{n=0}^\infty$
satisfies the weak orthogonal property
      \begin{align}\label{WeakOrtho1}
      \scal{\mathfrak{G}^{m,n} , \mathfrak{G}^{j,k}  }_{\omega}=
      0, \qquad \mbox{whenever }\quad n\ne k .
      \end{align}
Then the following identities hold
      \begin{align}
      &\scal{z \mathfrak{G}^{m,n},\mathfrak{G}^{j,k}}_{\omega} =0 ;
      \qquad\qquad\qquad ~ k\ne n  ~\mbox{ and }~ k\ne n-1
      \label{RecInt2}
      \\
      &\scal{z \mathfrak{G}^{m,n},\mathfrak{G}^{j,n}}_{\omega} =
      \scal{\mathfrak{G}^{m+1,n},\mathfrak{G}^{j,n}}_{\omega}
      \label{RecInt4}
      \\
      &\scal{z \mathfrak{G}^{m,n+1},\mathfrak{G}^{j,n}}_{\omega} =  (n+1)
      \scal{\mathfrak{G}^{m,n},\mathfrak{G}^{j,n}}_{\omega}
      \label{RecInt6}
      \end{align}
\end{proposition}

\begin{remark}
Note that one obtains by conjugation similar identities to
(\ref{RecInt2}), (\ref{RecInt4}) and (\ref{RecInt6}) with $\bz$
instead of $z$.
\end{remark}

\noindent {\sf Proof of Proposition \ref{PropId1}.} By multiplying both sides of the three-term
recursion relation (\ref{RecursionRelation2}),
      \begin{align*}
      \mathfrak{G}^{m+1,n}(z,\bz|\g)= z \mathfrak{G}^{m,n}(z,\bz|\g) - n
      \mathfrak{G}^{m,n-1}(z,\bz|\g), \label{RecRel2}
      \end{align*}
by $[\mathfrak{G}^{j,k}]^{*}$ and integrating over the whole $\C$
w.r.t. $\omega d\lambda$, we get
      \begin{equation}
      \scal{\mathfrak{G}^{m+1,n},\mathfrak{G}^{j,k}}_{\omega}= \scal{z
      \mathfrak{G}^{m,n},\mathfrak{G}^{j,k}}_{\omega} - n
      \scal{\mathfrak{G}^{m,n-1},\mathfrak{G}^{j,k}}_{\omega}.
      \label{RecInt1}
      \end{equation}
Next, using the weak orthogonality assumption (\ref{WeakOrtho1}),
together with (\ref{RecInt1}), we deduce easily the first one
(\ref{RecInt2}). Identities (\ref{RecInt4}) and (\ref{RecInt6}) are
obtained by taking respectively $k=n$ and $k=n-1$ in (\ref{RecInt1})
and applying again (\ref{WeakOrtho1}). \fin

Thus, making use of the recursion relation (\ref{RecursionRelation})
as well as of the above obtained identities, under the assumption
(\ref{WeakOrtho1}), one gets the following
      \begin{align}
      & \nu\scal{ \mathfrak{G}^{m,n}, \mathfrak{G}^{j+1,n}}_{\omega} +
      \frac{\bg}2\scal{ \mathfrak{G}^{m,n}, \mathfrak{G}^{j,n}}_{\omega}=
      m \scal{\mathfrak{G}^{m-1,n}, \mathfrak{G}^{j,n}}_{\omega}
      \label{RecInt1-5}
      \\ &
      \nu\scal{\bz \mathfrak{G}^{m,n}, \mathfrak{G}^{j,n+1}}_{\omega}=
      \scal{\mathfrak{G}^{m,n+1}, \mathfrak{G}^{j,n+1}}_{\omega}
      \label{RecInt1-6}
      \\ &
      \nu(n+1)\scal{\mathfrak{G}^{m,n}, \mathfrak{G}^{j,n}}_{\omega} =
      \scal{\mathfrak{G}^{m,n+1}, \mathfrak{G}^{j,n+1}}_{\omega}.
      \label{RecInt1-7}
      \end{align}
As consequence (together with the use of Proposition \ref{PropId1}),
we obtain useful identities for the norms. More precisely, we have
\begin{proposition}
Assume the weak orthogonality assumption (\ref{WeakOrtho1}) to be
satisfied by the set 
$\set{\mathfrak{G}^{m,n}:=\mathfrak{G}_{\nu}^{m,n}(\cdot,\cdot|\g)}_{n=0}^\infty$ for every fixed integer $m$.
Then we have
      \begin{align}
     i) & \quad \norm{\mathfrak{G}^{m,n+1}}^2_{\omega}= \nu^{n+1}(n+1)!
      \norm{\mathfrak{G}^{m,0}}^2_{\omega}. \label{RecInt1-10}
\\
      ii) & \quad \frac{\bg}2\scal{
      \mathfrak{G}^{m,n}, \mathfrak{G}^{m-1,n}}_{\omega} = \nu^{n}n! \big(
      m \norm{\mathfrak{G}^{m-1,0}}^2_{\omega} - \nu
      \norm{\mathfrak{G}^{m,0}}^2_{\omega}\big) . \label{RecInt1-12}
\\
iii) & \quad \mbox{The considered system is orthogonal if and only
if $\g=0$.} \nonumber
      \end{align}
\end{proposition}
\noindent {\sf Proof.} Putting $j=m$ in  (\ref{RecInt1-7}) infers
      \begin{align*}
      \nu(n+1)\norm{\mathfrak{G}^{m,n}}^2_{\omega} =
      \norm{\mathfrak{G}^{m,n+1}}^2_{\omega}.
      \end{align*}
Hence repeated application of the previous fact yields
      \begin{align*}
      \norm{\mathfrak{G}^{m,n+1}}^2_{\omega}= \nu^{l+1}(n+1)n \cdots (n+1-l)\norm{\mathfrak{G}^{m,n-l}}^2_{\omega}
      \end{align*}
for every positive integer $l\leq n$. In particular for $l=n$ we get
the asserted result (\ref{RecInt1-10}). While the result
(\ref{RecInt1-12}) in ii) holds by substitution of
(\ref{RecInt1-10}) in
      \begin{align}\label{RecInt1-8}
      \nu\norm{ \mathfrak{G}^{m,n}}^2_{\omega} + \frac{\bg}2\scal{
      \mathfrak{G}^{m,n}, \mathfrak{G}^{m-1,n}}_{\omega}= m
      \norm{\mathfrak{G}^{m-1,n}}^2_{\omega},
      \end{align}
that follows from (\ref{RecInt1-5}) for the specified value $j=m-1$.

\noindent For iii) suppose that
$\set{\mathfrak{G}^{m,n}}_{m,n=0}^\infty$ is orthogonal, i.e.,
$\scal{\mathfrak{G}^{m,n}  , \mathfrak{G}^{j,k} }_{\omega}= 0$
whenever $m\ne j$ or $n\ne k$. Then (\ref{RecInt1-5}) reduces
further for $j=m$ to $
\frac{\bg}2\norm{ \mathfrak{G}^{m,n}}^2_{\omega}=0,$ 
and therefore $\g=0$. The converse is obvious, indeed for $\g=0$,
the polynomials $\mathfrak{G}^{m,n}$ reduce to the complex Hermite
polynomials $H^{m,n}$ which are known to form an orthogonal system
w.r.t. the Gaussian measure $e^{-\nu |z|^2}d\lambda$. \fin

\begin{remark}
According to the classical fact that eigenfunctions associated to
different eigenvalues of a Hermitian operator are orthogonal, we
conclude from (\ref{EigenF}) that for every fixed positive integer
$m$ the set
$\set{\mathfrak{G}_{\nu}^{m,n}(\cdot,\cdot|\g)}_{n=0}^\infty$ is
orthogonal w.r.t. the weighting function $\omega(z)=e^{-\Re
e\scal{z,S_{\nu,\g}}}$. Hence, the weak orthogonality assumption
(\ref{WeakOrtho1}) is satisfied and therefore related identities
(\ref{WeakOrtho1})-(\ref{RecInt1-12}) hold.
\end{remark}

Now, let denote by $\norm{~}_{\nu,\g}$ the norm corresponding to
$\omega(z)=e^{-\Re e\scal{z,S_{\nu,\g}}}
$. We then assert

\begin{proposition}
The norm of the GCHP, $\mathfrak{G}^{m,n}$, w.r.t.
$\norm{~}_{\nu,\g}$ is given explicitly by
      \begin{align}\label{RecInt1-11}
      \norm{\mathfrak{G}^{m,n}}^2_{\nu,\g}= m!n!\pi \frac{\nu^{n}}{
      \nu^{m+1}} \mathbf{F}(m+1;1;\frac{|\g|^2}{4\nu}),
      \end{align}
where $\mathbf{F}(a;c;x)$ denotes the usual confluent hypergeometric
function.
\end{proposition}

\noindent {\sf Proof.}  The result in such proposition is a
particular case of the following Lemma, whose the proof can be
handled by straightforward computation. \fin

\begin{lemma} \label{LemExpNorm}
Let $C^\nu_{m,j}(\g)$ stands for
      \[
      C^\nu_{m,j}(\g):= \frac{(-1)^{m+j}\max(m,j)!\pi}{
      \nu^{\max(m,j)+1}2^{|m-j|}(|m-j|)!}\g^{\max(m,j)-j}{\bg}^{\max(m,j)-m}.
      \]
Then, we have
      \begin{align}
      \scal{\mathfrak{G}_{\nu}^{m,0} ,\mathfrak{G}_{\nu}^{j,0} }_{\nu,\g} &=
      \int_{\C} z^m {\bz}^j e^{-\nu|z|^2-\Re e\scal{z,\g}}d\lambda(z) \nonumber\\
      &= C^\nu_{m,j}(\g) \mathbf{F}(1+\max(m,j);1+|m-j|;\frac{|\g|^2}{4\nu}).\label{ExpComp}
      \end{align}
\end{lemma}

\begin{remark}
i) For the particular case of $\nu=1$ and $\g=0$, we recover from
(\ref{RecInt1-11}) the known fact that $\norm{H^{m,n}}^2= m!n!\pi $
for the complex Hermite polynomials \cite{IntInt06}.

ii) By combining  (\ref{RecInt1-11}) and (\ref{RecInt1-12}), we
conclude that
      \begin{align}
      \scal{
      \mathfrak{G}^{m,n}, \mathfrak{G}^{m-1,n}}_{\nu,\g} & =
      \frac{2m!n!\pi\nu^{n}}{ {\bg}\nu^{m}}\bigg( \mathbf{F}(m;1;\frac{|\g|^2}{4\nu})-
      \mathbf{F}(m+1;1;\frac{|\g|^2}{4\nu})\bigg). \label{RecInt1-13}
      \end{align}

iii) The explicit expression of $\scal{\mathfrak{G}_{\nu}^{m,0} ,
\mathfrak{G}_{\nu}^{j,0} }_{\nu,\g}$, 
given through (\ref{ExpComp}),  proves (again) that the family
$\set{\mathfrak{G}_{\nu}^{m,n}(\cdot,\cdot|\g)}_{m,n=0}^\infty$ is
not orthogonal w.r.t. $e^{-\Re e\scal{z,S_{\nu,\g}}}d\lambda$.
\end{remark}

\section{ Concluding remarks}

\subsection{ Additional identities for complex Hermite and Laguerre
polynomials}
\noindent We begin with the following for complex Hermite
polynomials
      \begin{align}\label{Identity2221}
      H^{m,n}(z,\nu\bz+\frac \bg 2)=
      \frac{(-1)^{(\mi)}n!}{\sqrt{\nu}^m}\sum_{j=0}^n
      \frac{\sqrt{\nu}^j}{j!}
      \frac{(\bg/2)^{n-j}}{(n-j)!}H^{m,j}(\sqrt{\nu} z,\sqrt{\nu}\bz)
      \end{align}
which follows from (\ref{Identity21}) and (\ref{LinkDattoli1}). It
yields in particular the following one
      \begin{align}
      H^{m,n}(z,\bz+1)= (-1)^{\mi}
      \sum_{j=0}^n\frac{n!}{j!(n-j)!}H^{m,j}(z,\bz).
      \label{Identity2222}
      \end{align}
An other identity that can we derive for Laguerre polynomials is
      \begin{align}
      L_{\mi}^{|m-n|}(z(\nu\bz +\frac{\bg}{2}))
      =  & \frac{(-1)^{(\mi)}n!}{\sqrt{\nu}^m(\mi)!}
      z^{(n-m)\curlyvee 0} (\nu\bz +\frac{\bg}{2})^{(m-n)\curlyvee 0}\times\nonumber\\
      &\times \sum_{j=0}^n
      \frac{\sqrt{\nu}^j}{j!}
      \frac{(\bg/2)^{n-j}}{(n-j)!}H^{m,j}(\sqrt{\nu}
      z,\sqrt{\nu}\bz). \label{Identity23}
      \end{align}
It follows from (\ref{LikeLaguerre}) combined with
(\ref{Identity21}), and gives rise in particular to
      \begin{equation}\label{Identity241}
      L_{n}^{m-n}(z\bz+z) = (-1)^{n}z^{n-m} \sum_{j=0}^n
      \frac{H^{m,j}(z,\bz)}{j!(n-j)!}
      \end{equation}
 for $m\geq n$ and
      \begin{equation}\label{Identity242}
      L_{m}^{n-m}(z\bz+z)   = (-1)^{m}
      \frac{n!}{m!}(\bz+1)^{m-n}    \sum_{j=0}^n
      \frac{H^{m,j}(z,\bz)}{j!(n-j)!}
      \end{equation}
 for $n\geq m$.\\
From (\ref{Identity241}) and (\ref{Identity242}) one can deduce also
the following
      \begin{equation}
      {n!} (\bz+1)^{m-n} \sum_{j=0}^n\frac{H^{m,j}(z,\bz)}{j!(n-j)!}
      =  {m!}z^{m-n}    \sum_{j=0}^m \frac{H^{n,j}(z,\bz)}{j!(m-j)!} ,
      \end{equation}
where we have assuming $m\geq n$.

\subsection{Illustration: Matrix representation}\label{MatrixRep}
Here we illustrate some obtained results making use of the matrix
representation of a polynomial
      \[
      P^{m,n}(z,\bz)=\sum_{j,k=0}^{m,n}p_{jk}z^j{\bz}^k
      \]
of degree $m$ in $z$ and degree $n$ in $\bz$. Thus the entries of
the matrix representing $P^{m,n}(z,\bz)$ are the coefficients
$p_{jk}$ (corresponding to the monomial $z^j{\bz}^k$) in such
expansion, i.e.,
      \[
      P^{m,n}(z,\bz)=
      \begin{array}{ccccccccccccccccc}
      1   \\  z  \\  \vdots \\   z^m
      \end{array}
      \stackrel{\begin{array}{cccccccccccc} \qquad ~ 1~~  & ~  \bz ~&
      \cdots\cdots& ~{\bz}^{n}\quad
      \end{array}}
      {\left(~\begin{array}{ccccccccccc}
      p_{00} &  p_{01} &  \cdots \cdots &  p_{0n}  \\
      p_{10} &  p_{11} &  \cdots \cdots &  p_{1n}  \\
      \vdots &  \vdots &                &  \vdots  \\
      p_{m0} &  p_{m1} &  \cdots \cdots &  p_{mn}
      \end{array}~ \right)}.
      \]
In the sequel we consider only the polynomials
$\mathfrak{G}_{\nu}^{m,n}(z,\bz|\g)$ with $m=n$, so that one deals
with  square matrices. In this case the polynomials
$\mathfrak{G}_{\nu}^{m,m}(z,\bz|\g)$ reduce further to
      \begin{align}\label{LikeHermitePoly2z}
      \mathfrak{G}_{\nu}^{m,m}(z,\bz|\g)
      = (-1)^{m} (m!)^2
      \sum_{l=0}^{m}\sum_{k=0}^{l} \frac{(-\nu)^l}{l!(m-l)!}
      \frac{\Big(\frac \bg {2\nu}z\Big)^{l-k}}{(l-k)!}
      \frac{|z|^{2k}}{k!},
      \end{align}
and one asserts that the monomials ${\bz}^j|z|^{2k}$, $j\ne 0$, do
not appear in the expansion. Hence
$\mathfrak{G}_{\nu}^{m,m}(z,\bz|\g)$ are represented by triangular
square matrices $\big[ g_{lk}(\bg) \big]_{l,k=0,1,\cdots, m}$, whose
entries $g_{lk}(\bg)$ are given by
      \[
      g_{lk}(\bg)=(-1)^m(m!)^2\left\{\begin{array}{ll}
      {\large\mbox{$\frac{(-1)^l\nu^k}{l!(m-l)!k!(l-k)!}$}}
      \Big(\frac{\bg}{2}\Big)^{l-k} &  \quad \mbox{if } ~ k\leq l \\
      0 &  \quad \mbox{if } ~ k> l \end{array}\right.
      \]
Furthermore, one can remark that the complex Hermite polynomials
$H^{m,k}(\sqrt\nu z,\sqrt\nu \bz)$, $k=m,m-1,\cdots,0$,  up to a
precise multiplicative constant $C_{m,k}(\bg)$, are respectively the
successive diagonals of the polynomial
$\mathfrak{G}_{\nu}^{m,m}(z,\bz|\g)$. In fact this observation is
contained in the formula (\ref{Identity21}). For illustration, we
consider the case where $\nu=1$ and $\g=2$ for which
(\ref{Identity21}) reads simply as
      \begin{equation}
      \label{Identity} \mathfrak{G}_{1}^{m,m}(z,\bz|2) = \sum_{k=0}^m
      \frac{m!}{k!(m-k)!}\cdot H^{m,k}(z,\bz) .
      \end{equation}
Added to $\mathfrak{G}_{1}^{0,0}(z,\bz|2)=\big(1\big)$ the first few
of these polynomials are given by
      \[
      \stackrel{\begin{array}{c}\large{\mbox{$\mathfrak{G}_{1}^{2,2}(z,\bz|2)$}}\\
      \shortparallel\end{array}} { \left(
      \begin{array}{cc} \fbox{-1} & 0 \\ ~~ 1 & \fbox{1} \end{array}\right)},
\qquad
      \stackrel{\begin{array}{c}\large{\mbox{$\mathfrak{G}_{1}^{2,2}(z,\bz|2)$}}\\
      \shortparallel\end{array}} { \left(
      \begin{array}{ccc} ~~\fbox{2} & 0 & 0\\  -4 & \fbox{-4} & 0 \\ ~~{1} & 2 &
      \fbox{1}
      \end{array}\right)}
\quad \mbox{and }
      \stackrel{\begin{array}{c}\large{\mbox{$\mathfrak{G}_{1}^{3,3}(z,\bz|2)$}}\\
      \shortparallel\end{array}} {\left(
      \begin{array}{cccc}
      \fbox{-6} & 0 & 0 & 0 \\
      {18} & \fbox{18} & 0 & 0 \\
      {-9} & -18 & \fbox{-9} & 0 \\
      {1} & {3} & {3} & \fbox{1}
      \end{array}\right)} .
      \]
Their analogues for the complex Hermite polynomials are
$H^{0,0}(z,\bz)=\big(1\big)$,
      \[
       \stackrel{\begin{array}{c}\large{\mbox{$H^{1,1}(z,\bz)$}}\\
      \shortparallel\end{array}}
      {\left(
      \begin{array}{cc} \fbox{-1} & 0 \\ ~~ 0 & \fbox{1} \end{array}\right)},
\quad
      \stackrel{\begin{array}{c}\large{\mbox{$H^{2,2}(z,\bz)$}}\\
      \shortparallel\end{array}}{\left(
      \begin{array}{ccc} ~~\fbox{2} & 0 & 0\\  0 & \fbox{-4} & 0 \\ 0 & 0 &
      \fbox{1}
      \end{array}\right)}
\quad \mbox{and }
      \stackrel{\begin{array}{c}\large{\mbox{$H^{3,3}(z,\bz)$}}\\ \shortparallel\end{array}}
      {\left(
      \begin{array}{cccc}
      \fbox{-6} & 0 & 0 & 0 \\
      0 & \fbox{18} & 0 & 0 \\
      0 & 0 & \fbox{-9} & 0 \\
      0 & 0 & 0 & \fbox{1}
      \end{array}\right)}.
      \]
Hence, the complex Hermite polynomial $H^{m,m}(z,\bz)$ appears as
the principal diagonal of its analogue
$\mathfrak{G}_{1}^{m,m}(z,\bz|2)$ (which also is the first column in
such matrix representation). Indeed,
      \begin{align*}
      \rm{Diag} \bigg[\mathfrak{G}_{1}^{m,m}(z,\bz|2)\bigg] & =
      (-1)^m(m!)^2\bigg(h_{kk}= \frac{(-1)^k}{(k!)^2(m-k)!} \bigg)_{k=0,1,\cdots, m} \\
      & = H^{m,m}(z,\bz).
      \end{align*}
This can be deduced also from (\ref{Identity}). Further, one notes,
when taking $m=3$ for example,  that added to $1\cdot
H^{3,3}(z,\bz)$ given above, we have
      \[
       \stackrel{\begin{array}{c}\large{\mbox{$3\cdot H^{3,2}(z,\bz)$}}\\ \shortparallel\end{array}}
      { \left(
      \begin{array}{cccc}
       0  &  0   & 0 & 0 \\
      18  &  0   & 0 & 0 \\
       0  & -18  & 0 & 0 \\
       0  &  0   & 3 & 0
      \end{array}\right)},
\qquad
       \stackrel{\begin{array}{c}\large{\mbox{$3\cdot H^{3,1}(z,\bz)$}}\\ \shortparallel\end{array}}
       {\left(
      \begin{array}{cccc}
       0  &  0   & 0 & 0 \\
       0  &  0   & 0 & 0 \\
      -9  &  0   & 0 & 0 \\
       0  &  3   & 0 & 0
      \end{array}\right)},
\qquad
       \stackrel{\begin{array}{c}\large{\mbox{$1\cdot H^{3,0}(z,\bz)$}}\\ \shortparallel\end{array}}
       {\left(
      \begin{array}{cccc}
       0  &  0   & 0 & 0 \\
       0  &  0   & 0 & 0 \\
       0  &  0   & 0 & 0 \\
       1  &  0   & 0 & 0
      \end{array}\right)}
      \]
and therefore the complex Hermite polynomials $ H^{3,k}(z,\bz)$,
$k=3,2,1,0$, are respectively the successive diagonals of the
polynomial $\mathfrak{G}_{1}^{3,3}(z,\bz|2)$.

An other fact to be signaled here is linked to the established fact
in Eqn. (b),
      \[
      \Pag \mathfrak{G}_{\nu}^{m,n}(z,\bz|\g)=\frac
      n2\mathfrak{G}_{\nu}^{m,n-1}(z,\bz|\g),
      \]
which can be used to obtain $\mathfrak{G}_{\nu}^{m,m-1}(z,\bz|\g)$
from $\mathfrak{G}_{\nu}^{m,m}(z,\bz|\g)$ (and in general
$\mathfrak{G}_{\nu}^{m,n-1}(z,\bz|\g)$ from
$\mathfrak{G}_{\nu}^{m,n}(z,\bz|\g)$).  It is more practicable than
Eqn.(d) when the matrix representation is used.
 Indeed, we have the following prescription: by dropping the last
column
 of $\mathfrak{G}_{\nu}^{m,n}(z,\bz|\g)$
and next differentiating all the entries of the obtained matrix
w.r.t. $\bg$, we get the matrix representing
$(m/2)\mathfrak{G}_{\nu}^{m,m-1}(z,\bz|\g)$.
 For illustration, we
take for example $m=4$. The sum of the complex Hermite polynomials
$H^{m,k}(\sqrt\nu z,\sqrt\nu \bz)$, $k=4,3,2,1,0$, according to
(\ref{Identity21}), infers the GCHP
$\mathfrak{G}_{\nu}^{4,4}(z,\bz|\g)$ given by
      \[
       \left(~
      \begin{array}{ccccccccccc}
      \mathbf{24}                   &~&  0                               &~&  0                               &~&  0                  &~&   0    \\
      -96\Big(\frac{\bg}2\Big)      & &  \mathbf{-96\nu}                 & &  0                               & &  0                  & &  0     \\
      ~ 72{\Big(\frac{\bg}2\Big)}^2 & &  144\nu\Big(\frac{\bg}2\Big)     & &  \mathbf{72\nu^2}                & &  0                  & &  0     \\
      -16(\bg/2)^3                  & &  -48\nu{\Big(\frac{\bg}2\Big)}^2 & &  -48\nu^2\Big(\frac{\bg}2\Big)   & &  \mathbf{-16 \nu^3} & &  0     \\
      {\Big(\frac{\bg}2\Big)}^4     & &  4\nu{\Big(\frac{\bg}2\Big)}^3   & &  6\nu^2{\Big(\frac{\bg}2\Big)}^2 & &  4 \nu^3\Big(\frac{\bg}2\Big)
                                                                                                                                      & & \mathbf{\nu^4}
      \end{array}~ \right)
      \]
and therefore  (from the fact $\Pag
\mathfrak{G}_{\nu}^{4,4}(z,\bz|\g)=\frac 4 2 \cdot
\mathfrak{G}_{\nu}^{4,3}(z,\bz|\g)$) we get
      \[
      \mathfrak{G}_{\nu}^{4,3}(z,\bz|\g)=\left(~
      \begin{array}{ccccccccccc}
      0                             &~&  0                             &~&  0                           &~&  0     &     \\
      -24                           & &  0                             & &  0                           & &  0     &     \\
      ~ 36\Big(\frac{\bg}2\Big)     & &  36\nu                         & &  0                           & &  0     &     \\
      -12{\Big(\frac{\bg}2\Big)}^2  & &  -24\nu\Big(\frac{\bg}2\Big)   & &  -12\nu^2                    & &  0     &     \\
      {\Big(\frac{\bg}2\Big)}^3     & &  3\nu{\Big(\frac{\bg}2\Big)}^2 & &  3\nu^2\Big(\frac{\bg}2\Big) & &  \nu^3 &
      \end{array}~ \right).
      \]

\newpage

      \noindent \centerline{\bf Acknowledgments.}
       {The author is indebted to Professor S. Thangavelu
      for the interest that he yields for this work
      and to Professor M. Lassalle for his suggestions which led the
      paper its final form. Special thanks are addressed to
      Professor A. Intissar for his encouragements. The author
      would like to thank also the anonymous referee
      for bringing W\"{u}nsche papers to his attention.}

\end{document}